\documentclass{article}

    \PassOptionsToPackage{numbers, compress}{natbib}

\usepackage[final]{neurips_2024}

\usepackage{amssymb, amsmath}
\usepackage{stmaryrd}

\usepackage{hyperref}

\usepackage{bm}
\usepackage[textsize=tiny, disable]{todonotes}
\usepackage{soul}

\usepackage{tikz}
\usetikzlibrary{positioning}
\usepackage{pgfplots}
\pgfplotsset{compat=1.10}
\usepgfplotslibrary{groupplots,fillbetween}
\usetikzlibrary{spy, backgrounds}
\usetikzlibrary{shapes.geometric, arrows.meta}
\usepackage{diagbox}
\usepackage{hhline}

\newcommand{\rR}{\mathbb{R}}

\newcommand{\rC}{\mathbb{C}}

\newcommand{\bB}{\mathbf{B}}
\newcommand{\bC}{\mathbf{C}}

\newcommand{\bL}{\mathbf{L}}

\newcommand{\bK}{\mathbf{K}}

\newcommand{\bn}{\mathbf{n}}
\newcommand{\bq}{\mathbf{q}}
\newcommand{\bx}{\mathbf{x}}
\newcommand{\bu}{\mathbf{u}}

\newcommand{\tbq}{\tilde{\bq}}
\newcommand{\tbu}{\tilde{\bu}}
\newcommand{\tp}{\tilde{p}}

\newcommand{\bp}{\mathbf{p}}
\newcommand{\f}{\mathbf{f}}
\newcommand{\bPhi}{\mathbf{\Phi}}
\newcommand{\bphi}{\bm{\phi}}

\newcommand{\hC}{{\hat{C}}}

\newcommand{\hu}{{\hat{u}}}

\newcommand{\hbu}{{\hat{\bu}}}

\newcommand{\hbp}{{\hat{\bp}}}

\newcommand{\hbB}{{\hat{\bB}}}
\newcommand{\hbC}{{\hat{\bC}}}

\newcommand{\hbK}{{\hat{\bK}}}
\newcommand{\hbL}{{\hat{\bL}}}

\newcommand{\bOmega}{\overline{\Omega}}

\def\inprod#1#2{\left\langle #1, #2 \right\rangle}

\def\jump#1{\left\llbracket #1 \right\rrbracket}
\def\avg#1{\left\{\!\!\left\{ #1 \right\}\!\!\right\}}

\pgfplotsset{select coords between index/.style 2 args={
    x filter/.code={
        \ifnum\coordindex<#1\fi
        \ifnum\coordindex>#2\fi
    }
}}

\usepackage[utf8]{inputenc} 
\usepackage[T1]{fontenc}    
\usepackage{hyperref}       
\usepackage{url}            
\usepackage{booktabs}       
\usepackage{amsfonts}       
\usepackage{nicefrac}       
\usepackage{microtype}      
\usepackage{xcolor}         

\title{Scalable physics-guided data-driven component model reduction for steady Navier-Stokes flow}

\author{%
Seung Whan Chung$^1$
\quad
Youngsoo Choi$^1$
\quad
Pratanu Roy$^2$
\\
\textbf{Thomas Roy}$^3$
\quad
\textbf{Tiras Lin}$^3$
\quad
\textbf{Du T. Nguyen}$^4$
\\
\textbf{Christopher Hahn}$^4$
\quad
\textbf{Eric B. Duoss}$^5$
\quad
\textbf{Sarah E. Baker}$^4$
\\
$^1$Center for Applied Scientific Computing
\qquad
$^2$Atmospheric, Earth and Energy Division
\\
$^3$Computational Engineering Division
\qquad
$^4$Material Science Division
\\
$^5$Material Engineering Division
\\
Lawrence Livermore National Laboratory, Livermore, CA 94550 \\
\texttt{\{chung28,choi15,roy23,roy27,lin46,}
\\
\texttt{hahn31,duoss1,baker74\}@llnl.gov}
\\
\texttt{dunguyen@gmail.com}
}

\newif\ifcommentcoloring
\commentcoloringfalse

\newcommand{\KC}[1]{\ifcommentcoloring \textcolor{brown!70!black}{#1}\else #1\fi}

\begin{document}

\maketitle

\begin{abstract}
Computational physics simulation can be a powerful tool to accelerate
industry deployment of new scientific technologies.
However, it must address the challenge of
computationally tractable, moderately accurate prediction at large industry scales,
and training a model without data at such large scales.
A recently proposed component reduced order modeling (CROM) tackles this challenge by combining reduced order modeling (ROM)
with discontinuous Galerkin domain decomposition (DG-DD).
While it can build a component ROM at small scales that can be assembled into a large scale system,
its application is limited to linear physics equations.
In this work, we extend CROM to nonlinear steady Navier-Stokes flow equation.
Nonlinear advection term is evaluated via tensorial approach or empirical quadrature procedure.
Application to flow past an array of objects at moderate Reynolds number demonstrates
$\sim23.7$ times faster solutions with a relative error of  $\sim 2.3\%$, even at scales
$256$ times larger than the original problem.
\end{abstract}

\section{Introduction}

Industry deployment of a novel scientific technology often involves scaling up process,
which demonstrates performance of a lab-scale proven method at industry scale.
Conventionally, the scaling up process is performed through physical pilot plants at intermediate scales,
though they are costly and time-consuming to design, construct and operate.
Computational simulations can augment and accelerate design process and prediction in this deployment procedure.
However, even pilot scales are often order-of-magnitude larger than lab scale,
which is computationally intractable with traditional numerical methods.
Sub-grid scale approximations such as volume-averaging or closure model for large-eddy simulations can compromise the accuracy significantly~\cite{song2004prediction}.\todo{KC: maybe I should be careful about this sentence.}
Meanwhile, extremely large scale application also challenges use of recent data-driven methods,
in that there is no available data at such large scale and the prediction must be extrapolation in scale.
\par
The recently proposed component reduced order modeling (CROM)~\cite{chung2024train, chung2023scalable} tackles this challenge by combining projection-based reduced order modeling (PROM)
with discontinuous Galerkin domain decomposition (DG-DD).
Proper orthogonal decomposition (POD)~\cite{berkooz1993proper} identifies a low-dimensional linear subspace
that can effectively represent the physics solutions based on small scale sample snapshot data.
PROM projects the physics governing equation onto the linear subspace,
thereby achieving both robust accuracy and cheap computation time.
Small scale unit reduced-order models (ROMs) are then assembled into a large scale ROM system,
where the interface condition is handled via discontinuous Galerkin penalty terms.
Since linear subspace identification and the reduced order modeling can be performed only at the small unit scales,
CROM can achieve robust extrapolation in scale without data at large scale.
\par
CROM has been successfully demonstrated for several applications such as
Poisson equation, Stokes flow, advection-diffusion equation~\cite{chung2024train, chung2023scalable}, and linear elasticity~\cite{mcbane2021component, mcbane2022stress}.
However, all of these applications have been limited to linear systems.
In this work, we extend CROM to nonlinear equations, particularly steady incompressible Navier-Stokes equation.
Naive projection of nonlinear terms onto linear subspace would not gain any speed-up,
\KC{requiring an efficient approximation technique.}
We address this issue with two different approaches.
First, exploiting the fact that the advection is quadratic in terms of velocity,
we can pre-compute a 3rd-order tensor ROM operator for advection~\cite{Lassila2014}.
Second, we can also employ empirical quadrature procedure (EQP)~\cite{chapman2017accelerated}
to evaluate advection term only at the selected sample grid points,
which are obtained from a minimization problem with respect to sample data.
Furthermore, the incompressibility of the physics necessitates the linear subspaces to satisfy
the associated inf-sup condition~\cite{Babuvska1971,Brezzi1974,Ladyzhenskaya1963,Taylor1973}.
This can be addressed by augmenting the velocity bases with compressible components from gradients of pressure POD modes~\cite{ballarin2015supremizer}.
\par
The rest of the paper is organized as follows. In Section~\ref{sec:formulation}, we provide a concise overview of the proposed
component model reduction approach with the specific example of steady incompressible Navier-Stokes equation.
Following that, in Section~\ref{sec:result}, we demonstrate of the proposed method to a scaled-up prediction of flow past an array of objects at moderate Reynolds number.

\vspace{-10pt}
\section{Formulation}\label{sec:formulation}
\vspace{-10pt}

We consider the global-scale domain $\Omega\subset\rR^d$ decomposed into $M$ subdomains $\Omega_m$, i.e. $\Omega = \bigcup_{m=1}^M\Omega_m$.
All subdomains can be categorized into a few reference domains $\rC\equiv\left\{\bOmega_1, \bOmega_2, \ldots\right\}$.
Steady incompressible Navier-Stokes equation
for each subdomain velocity $\tbu_m \in H_1(\Omega_m)^d$ and pressure $\tp_m\in H_1(\Omega_m)$ writes
\begin{subequations}\label{eq:dd-gov}
    \begin{equation}
        -\nu\nabla^2\tbu_m + \nabla\tp_m + \tbu_m\cdot\nabla\tbu_m = \f_m
    \end{equation}
    \begin{equation}\label{eq:dd-gov-div}
        \nabla\cdot\tbu_m = 0,
    \end{equation}
    with non-dimensional viscosity $\nu=1/\mathrm{Re}$ as the inverse of Reynolds number.
    The interface $\Gamma_{m,n}\equiv\partial\Omega_m\cap\partial\Omega_n$ is constrained by the continuity and smoothness condition,
    \begin{equation}
        \jump{\tbu} = \jump{\tp} = 0 \qquad \text{on } \Gamma_{m,n}
    \end{equation}
    \begin{equation}
        \avg{\bn\cdot\nabla\tbu} = \avg{\bn\cdot\nabla\tp} = 0 \qquad \text{on } \Gamma_{m,n},
    \end{equation}
    with $\jump{\tbq} \equiv \tbq_m - \tbq_n$ and $\avg{\bn\cdot\nabla\tbq} \equiv \frac{1}{2}\left( \bn_m\cdot\nabla\tbq_m + \bn_n\cdot\nabla\tbq_n \right)$.
\end{subequations}
$\bn_m$ is the outward normal vector of the subdomain $\Omega_m$ and $\bn_m = -\bn_n$ on $\Gamma_{m,n}$.
\par
DG-DD seeks an approximate solution $(\bu, \bp)=\{(\bu_m, \bp_m)\}$ that satisfies the discretization of (\ref{eq:dd-gov}),
\begin{subequations}\label{eq:dd-gov-dg}
    \begin{equation}
        \bK_m\bu_m + \bB_m^{\top}\bp_m + \bC_m[\bu_m]
        + \sum\limits_{\Gamma_{m,n}\ne\emptyset}
        \left\{
        \begin{pmatrix}
            \bK_{mm} & \bK_{mn} \\
        \end{pmatrix}
        \begin{pmatrix}
            \bu_m \\ \bu_n
        \end{pmatrix}
        +
        \begin{pmatrix}
            \bB_{mn}^{\top} & \bB_{mn}^{\top} \\
        \end{pmatrix}
        \begin{pmatrix}
            \bp_m \\ \bp_n
        \end{pmatrix}
        \right\}
        = \bL_m,
    \end{equation}
    \begin{equation}
        \bB_m\bu_m
        + \sum\limits_{\Gamma_{m,n}\ne\emptyset}
        \begin{pmatrix}
            \bB_{mn} & \bB_{mn} \\
        \end{pmatrix}
        \begin{pmatrix}
            \bu_m \\ \bu_n
        \end{pmatrix}
        = 0,
    \end{equation}
    for $\forall m$.
\end{subequations}
$\bK_m$, $\bB_m$ and $\bC_m$ correspond to viscous flux, velocity divergence and nonlinear advection operator in the physics equation (\ref{eq:dd-gov}a-b), respectively.
The summations over interfaces $\Gamma_{m,n}$ weakly enforces the interface condition (\ref{eq:dd-gov}c-d).
This discretized physics equation provides the base for component ROM.
\par
We approximate the component-level solution $(\bu_r, \bp_r)$ on a low-dimensional linear subspace,
\begin{equation}\label{eq:pod}
    \bu_r \approx \bPhi_{u,r}{\hbu}_r
    \qquad\qquad
    \bp_r \approx \bPhi_{p,r}{\hbp}_r,
\end{equation}
where the reduced solution $({\hbu}_r, \hbp_r)$ are the coefficients of the column vectors of the basis $\bPhi_{u,r}$ and $\bPhi_{u,r}$.
The bases are identified from sample solutions of (\ref{eq:dd-gov-dg}) on the reference domains $\bOmega_r$
via POD~\cite{Chatterjee2000,Liang2002}.
\KC{
The basis size is determined so that the sampled solutions may be represented with the basis at a desired accuracy.
This accuracy can be evaluated by the ratio between the sum of the given basis vectors' singular values and the total sum of all singular values,
\begin{equation}\label{eq:pod-eps}
    \epsilon_R = 1 - \frac{\sum_s^{R}\sigma_s}{\sum_s^{S}\sigma_s}.
\end{equation}
This is equivalent to the relative representation error of the linear subspace spanned by the given basis over the sampled component-level solutions. In this study, both velocity and pressure POD basis sizes are chosen to be $40$, with $\sim3\%$ relative error over the sampled solutions.
}
Due to the incompressible nature of solution for (\ref{eq:dd-gov-dg}),
the velocity basis $\bPhi_{u,r}$ also remains divergence-free, resulting in spurious pressure modes.
To avoid this, the velocity basis is further augmented with the gradients of pressure POD modes as supremizer~\cite{ballarin2015supremizer},
\begin{equation}\label{eq:supremizer}
    \bPhi_{u,r} =
    GS\,
    \begin{pmatrix}
        \bPhi_{u,r}^{POD}
        &
        \bB_r\bPhi_{p,r}^{POD}
    \end{pmatrix},
\end{equation}
where $GS(\cdot)$ is modified Gram-Schmidt orthonormalization.
\par
With the linear subspace approximation (\ref{eq:pod}),
The physics equation (\ref{eq:dd-gov-dg}) is then projected onto the column space of $\bPhi_{u,r}$ and $\bPhi_{p,r}$,
\begin{subequations}\label{eq:dd-rom}
    \begin{equation}
        \hbK_m\bu_m + \hbB_m^{\top}\bp_m + \hbC_m[\bu_m]
        + \sum\limits_{\Gamma_{m,n}\ne\emptyset}
        \left\{
        \begin{pmatrix}
            \hbK_{mm} & \hbK_{mn} \\
        \end{pmatrix}
        \begin{pmatrix}
            \hbu_m \\ \hbu_n
        \end{pmatrix}
        +
        \begin{pmatrix}
            \hbB_{mn}^{\top} & \hbB_{mn}^{\top} \\
        \end{pmatrix}
        \begin{pmatrix}
            \hbp_m \\ \hbp_n
        \end{pmatrix}
        \right\}
        = \hbL_m,
    \end{equation}
    \begin{equation}
        \hbB_m\hbu_m
        + \sum\limits_{\Gamma_{m,n}\ne\emptyset}
        \begin{pmatrix}
            \hbB_{mn} & \hbB_{mn} \\
        \end{pmatrix}
        \begin{pmatrix}
            \hbu_m \\ \hbu_n
        \end{pmatrix}
        = 0,
    \end{equation}
\end{subequations}
for $\forall m$, with ROM operators
$\hat{\bK}_m=\bPhi_{u,m}^{\top}\bK_m\bPhi_{u,m}$, $\hat{\bK}_{ij} = \bPhi_{u,i}^{\top}\bK_{ij}\bPhi_{u,j}$,
$\hat{\bB}_m=\bPhi_{p,m}^{\top}\bB_m\bPhi_{u,m}$, $\hat{\bB}_{ij} = \bPhi_{p,i}^{\top}\bB_{ij}\bPhi_{u,j}$ and
$\hat{\bL}_m = \bPhi_{u,m}^{\top}\bL_m$.
The operators in (\ref{eq:dd-rom}) are the building blocks of the global ROM. The nonlinear ROM operator $\bC$ is described subsequently.
\par
Unlike linear ROM operators $\bK$ and $\bB$, the nonlinear operator in general cannot can be pre-computed as a reduced matrix.
In this work, the nonlinear ROM operator $\hbC$ is evaluated in two different approaches.
First, exploiting the fact that the advection term is quadratic with respect to $\bu$, we pre-compute a 3rd-order tensor operator,
\begin{equation}\label{eq:tensor}
        \hC_m[\hbu_m]_i
        = \sum_{j,k}\hbC_{ijk}\hu_{m,j}\hu_{m,k}
        \equiv \sum_{j,k}\inprod{{}_i\bphi_{u,m}}{{}_j\bphi_{u,m}\cdot\nabla{}_k\bphi_{u,m}}_{\Omega_m}\hu_{m,j}\hu_{m,k},
\end{equation}
for $\forall i\in[1,\dim(\hbu_m)]$,
where ${}_i\bphi_{u,m}$ is the $i$-th column vector of $\bPhi_{u,m}$ and $\inprod{\cdot}{\cdot}_{\Omega_m}$ is the inner product over the subdomain $\Omega_m$.
While the complexity of (\ref{eq:tensor}) scales faster than linear ROM operators,
we can still expect a significant speed-up if a moderate size of basis is used.
\par
An alternative is the empirical quadrature procedure (EQP) where the nonlinear term is evaluated at sampled grid points,
\begin{equation}\label{eq:eqp}
    \hC_m[\hbu_m]_i = \sum_q^{N_q} {}_i\bphi_{u,m}^{\top}(\bx_q)\;w_q\;[\overline{\bu}_m(\bx_q)\cdot\nabla\overline{\bu}_m(\bx_q)]
    \qquad \forall i\in[1,\dim(\hbu_m)],
\end{equation}
with $\overline{\bu}_m = \bPhi_{u,m}\hbu_m$.
Note that the summand in (\ref{eq:eqp}) is evaluated only at the EQP points $\bx_q$,
thus we can still expect similar speed-up with tensorial approach, as long as $N_q\sim \dim(\hbu_m)$.
Furthermore, this EQP approach is applicable to general nonlinear equations.
The EQP points $\bx_q$ and weight $w_q$ are calibrated with respect to velocity basis and sample solutions on the reference domain
via non-negative least-squares method~\cite{chapman2017accelerated}.
The number of points $N_q$ is controlled by error threshold of the quadrature,
which is set to $1\%$ in this study.

\vspace{-10pt}
\section{Results}\label{sec:result}
\vspace{-10pt}
We demonstrate CROM for steady Navier-Stokes flow
on the flow past array of objects used in Chung \textit{et al.}~\cite{chung2024train}.
Five different reference domains $\rC =\{\bOmega_1, \ldots, \bOmega_5\}$ are
considered as components for building up the global-scale system.
All reference domains lie within a unit square $\bOmega_r \subset [0, 1]^2$
with an obstacle within them: circle, square, triangle, star, and none (empty).
To obtain the POD bases, sample snapshots are generated on $2000$ 2-by-2-component domains with four randomly chosen subdomains from $\rC$.
The FOM (\ref{eq:dd-gov-dg}) is solved for $\nu=0.04$ ($\mathrm{Re}=25$) with the inflow velocity randomly chosen per each sample.
\KC{For flows past blunt bodies, steady flow physically exists only up to $\mathrm{Re}\lesssim40$ based on the object length scale~\cite{lienhard1966synopsis}.}
For the details of the reference domain meshes and sampling procedures,
we refer readers to Chung \textit{et al.}~\cite{chung2024train}.
The details of method implementation and instruction for the main experimental results can be found in the open-source code \href{https://github.com/LLNL/scaleupROM}{\texttt{scaleupROM}} (MIT license).
All numerical experiments are performed on an Intel Sapphire Rapids 2GHz processor with 256GB memory.
\par
In order to obtain the ROM bases (\ref{eq:pod}), POD is performed over velocity and pressure snapshots of each reference domain.
In this demonstration, 40 POD modes are chosen for both velocity and pressure,
\KC{which can represent overall snapshots with $\sim3\%$ relative error.}
The velocity POD bases are further augmented with the gradient of pressure POD modes per (\ref{eq:supremizer},
having additional 40 basis vectors.
\par
\begin{figure}[tbph]
  \centering
    \includegraphics[height=0.55\textwidth]{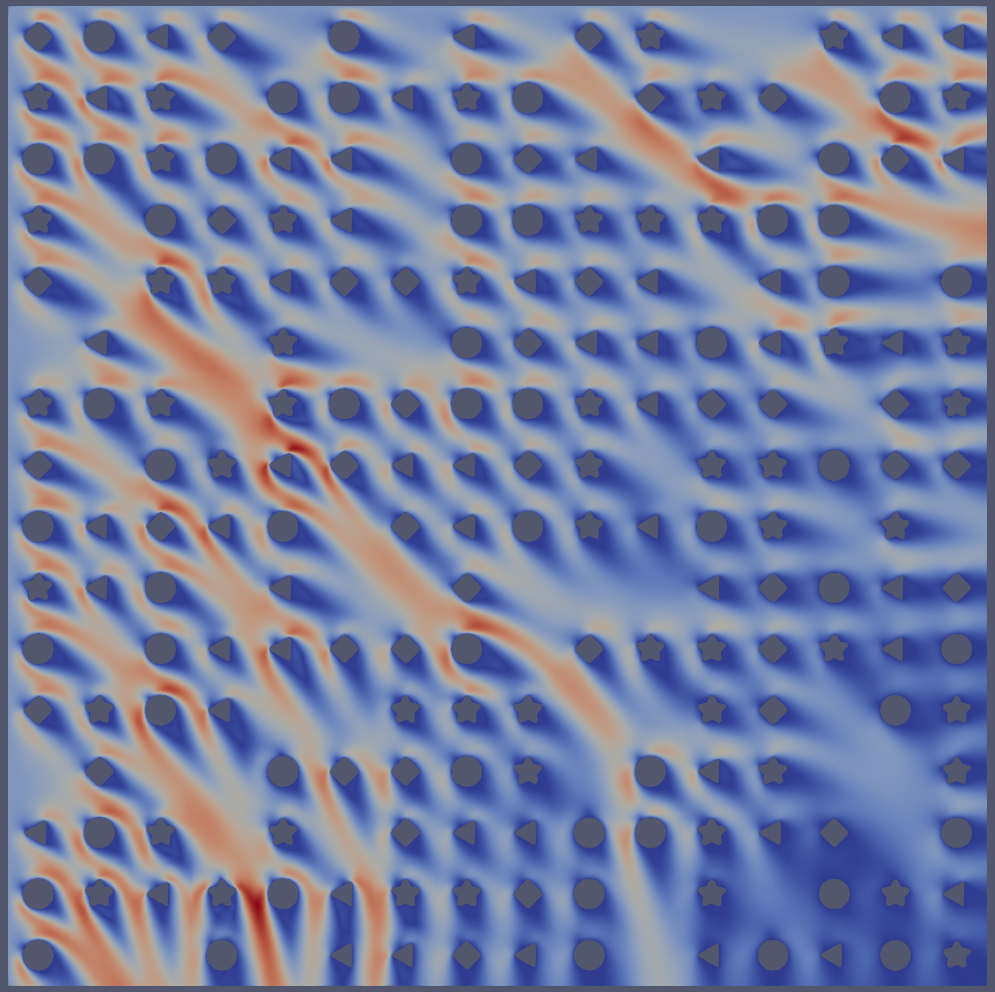}
    \includegraphics[height=0.55\textwidth]{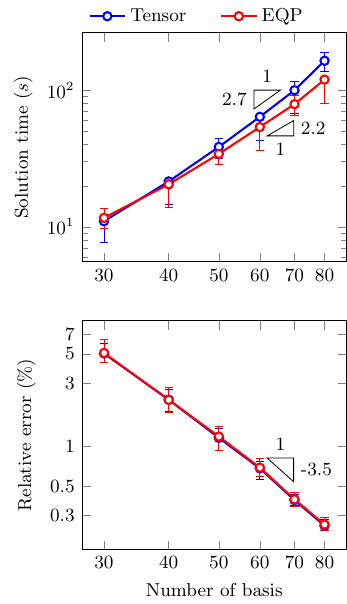}
    \caption{Scaled-up prediction of steady Navier-Stokes flow at $\nu=0.04$ over a $16\times16$ array of 5 different random objects:
    (left) \KC{flow-speed prediction of the proposed CROM};
    (top right) computation time with either (\ref{eq:tensor}) or (\ref{eq:eqp}), with respect to number of basis vectors; and
    (bottom right) relative error compared to the FOM solution, with respect to number of basis vectors.
    The error bar indicates $95\%$ confidence interval over 100 sample cases.
    }
    \label{fig:result}
\end{figure}
Using ROM (\ref{eq:dd-rom}),
scaled-up predictions are performed on $100$ test cases of $16\times16$ arrays of objects with random inflow velocity $\bu_{in}\in U[-1,1]^2$.
Figure~\ref{fig:result} shows an example scaled-up prediction for a $16\times16$ array of objects.
Given total $120$ basis vectors, ROM achieved about $23.7\times$ speed-up while maintaining $\sim2.3\%$ relative error compared to the FOM solution.
It is worth emphasizing that the ROM operators in (\ref{eq:dd-rom}) are built only from $2\times2$ domains.
On larger domains, the flow tends to accumulate on empty subdomains, which cannot be observed from the sampled snapshots.
However, ROM was able to robustly extrapolate in scale based on its underlying physics governing equation.
\par
Right subfigures of Figure~\ref{fig:result} compare two different approaches to evaluate the nonlinear advection.
With 30 basis vectors, both tensorial approach and EQP method takes similar computation time.
However, the computation time for EQP scales slightly better with number of basis vectors compared to tensorial approach.
Such speed-up did not come with a compromise in its accuracy almost at all, showing EQP's superiority over the tensorial approach.
Overall, for both approaches, the relative error scales much faster than the computation time, showing the effectiveness of ROM.
\par
\begin{table}
\begin{tabular}{|l||c|c|c|c|c|}
\hline
\diagbox{$R_{p}$}{$Z_p$} & 20 & 30 & 40 & 50 & 60 \\
\hline
20 & (2.5, 1.7) & (2.5, 3.5) & (2.4, 3.8) & (2.4, 3.6) & (2.4, 3.3) \\
\hline
30 & (2.8, 61.5) & (2.5, 1.3) & (2.5, 2.7) & (2.4, 2.6) & (2.4, 2.1) \\
\hline
40 & (7.5, $1.1\times10^5$) & (3.0, 50.9) & (2.6, 1.0) & (2.5, 1.6) & (2.4, 1.5) \\
\hline
50 & (28, $9.3\times10^6$) & (7.7, $6.8\times10^4$) & (3.0, 28.4) & (2.6, 1.1) & (2.5, 1.8) \\
\hline
60 & (N/A, N/A) & (26.2, $5.8\times10^6$) & (7.8, $4.7\times10^4$) & (3.1, 34.7) & (2.7, 1.2) \\
\hline
\end{tabular}
\caption{
\KC{Relative error of (flow velocity, pressure) in percentage,
depending on pressure basis size $R_p$ and supremizer size $Z_p$.
The predictions are made on a $16\times16$ array with velocity POD basis size of 40 (before augmentation).
N/A indicates that the numerical solution is not converged.}}
\label{table:ablation}
\end{table}
\KC{
A numerical experiment is further performed to investigate the impact of augmenting the velocity basis (\ref{eq:supremizer}).
Table~\ref{table:ablation} shows the relative errors of ROM predictions for a $16\times16$ array with different pressure basis sizes $R_p=\dim(\bPhi_{p,r})$ and supremizer sizes $Z_p=\dim(\mathbf{B}_r\bPhi_{p,r}^{POD})$.
For all cases of $R_p > Z_p$, the accuracy for the pressure degrades quickly with $Z_p$ due to spurious pressure modes.
This also impacts the accuracy for the flow velocity as well.
For all $R_p$, the pressure error is lowest at $R_p=Z_p$,
and plateaus for $R_p < Z_p$.
The plateau of the pressure error gradually decreases with $R_p$,
as more pressure basis vectors resolve the solution.
As long as more supremizers are used than pressure basis vectors, i.e. $R_p \le Z_p$,
the error for the flow velocity is maintained consistently,
as the same velocity POD basis is used over all cases.
This result strongly shows the role of the supremizer for stabilizing ROM pressure predictions.
}

\vspace{-10pt}
\section{Conclusion}\label{sec:conclusion}
\vspace{-10pt}

In order to accelerate a scaling-up process,
computational simulation must be able to reliably extrapolate in scale only from small-scale data.
CROM realizes this by combining ROM with DG-DD,
though it has been limited only to linear physics equations.
In this work we extend CROM to steady, nonlinear Navier-Stokes equation.
Nonlinear advection term is evaluated by either tensorial approach or EQP.
In order to avoid spurious pressure modes coming from divergence-free velocity bases,
the velocity bases are augmented with the gradients of pressure POD modes.
The proposed method is demonstrated on a porous media flow problem,
where 5 different unit components are combined into a large array.
The resulting ROM accelerates the solving time by a factor of $\sim23.7$ only with $\sim2.3\%$ relative error,
for a global domain $256$ times larger than the components.
While EQP shows its superiority over tensorial approach for both computational time and accuracy,
both are shown to be effective for reliable ROM prediction.
\par
Though only demonstrated on steady Navier-Stokes equation,
the proposed method is applicable to general nonlinear physics equation.
For examples, an advection-diffusion-reaction equation can be added in order to solve a mass transfer problem.
While tensorial approach is limited to polynomially nonlinear terms,
EQP is in general applicable to any type of nonlinear terms.
\par
In this work, the interface penalty term did not involve any nonlinear terms,
though it is not a general limitation of this proposed method.
In principle, tensorial approach or EQP method is readily extensible for interface penalty terms.
Such extension would be valuable for unsteady hyperbolic conservation laws,
further elevating the applicability of CROM.

\begin{ack}
This work was performed under the auspices of the U.S. Department of Energy
by Lawrence Livermore National Laboratory under contract DE-AC52-07NA27344
and was supported by Laboratory Directed Research and Development funding under project 22-SI-006.
LLNL-PROC-868977.

\end{ack}



\bibliographystyle{elsarticle-num} 
\bibliography{references}









\end{document}